\documentclass{amsproc}
\usepackage{amsmath,amssymb,amsfonts}
\usepackage{eucal}

\newtheorem{theorem}{Theorem}[section]
\newtheorem{lemma}[theorem]{Lemma}
\newtheorem{prop}[theorem]{Proposition}

\theoremstyle{definition}
\newtheorem{definition}[theorem]{Definition}

\newtheorem*{thank}{Acknowledgments}
\newtheorem*{note}{Note}

\theoremstyle{remark}

\numberwithin{equation}{section}

\input xy
\xyoption{all}

\newcommand{\Deltaop}{{\bf \Delta}^{op}}
\newcommand{\IDeltaop}{{\bf I \Delta}^{op}}
\newcommand{\colim}{\text{colim}}
\newcommand{\hocolim}{\text{hocolim}}
\newcommand{\holim}{\text{holim}}
\newcommand{\nerve}{\text{nerve}}
\newcommand{\Hom}{\text{Hom}}
\newcommand{\Map}{\text{Map}}
\newcommand{\SSets}{\mathcal{SS}ets}
\newcommand{\SSetst}{\mathcal{SS}ets^{\mathcal T}}

\newcommand{\SSetstm}{\mathcal{SS}ets^{\mathcal T_M}}
\newcommand{\LSSetst}{\mathcal {LSS}ets^{\mathcal T}}
\newcommand{\LSSetstm}{\mathcal {LSS}ets^{\mathcal T_M}}
\newcommand{\Tm}{\mathcal T_M}
\newcommand{\Sets}{\mathcal Sets}
\newcommand{\Algt}{\mathcal Alg^\mathcal T}
\newcommand{\Algtm}{\mathcal Alg^{\mathcal T_M}}
\newcommand{\SSpo}{\mathcal{SS} p_\mathcal O}
\newcommand{\SSpof}{\mathcal{SS} p_{\mathcal O,f}}
\newcommand{\SSpoc}{\mathcal{SS} p_{\mathcal O,c}}
\newcommand{\LSSpof}{\mathcal{LSS}p_{\mathcal O,f}}
\newcommand{\LSSpoc}{\mathcal{LSS}p_{\mathcal O,c}}
\newcommand{\SSp}{\mathcal{SS}p}
\newcommand{\LSSp}{\mathcal{LSS}p}
\newcommand{\alphau}{{\underline \alpha}}
\newcommand{\xu}{{\underline x}}
\newcommand{\Tocat}{\mathcal T_{\mathcal {OC}at}}
\newcommand{\SSetsTocat}{\mathcal {SS}ets^{\mathcal
T_{\mathcal {OC}at}}}
\newcommand{\LSSetsTocat}{\mathcal {LSS}ets^{\mathcal
T_{\mathcal {OC}at}}}
\newcommand{\SSetso}{\mathcal {SS}ets^{{\bf \Delta}^{op}_\mathcal
O}}

\newcommand{\Secat}{\mathcal Se \mathcal Cat}

\begin{document}

\title{Simplicial monoids and Segal categories}

\author{Julia E. Bergner}
\address{Kansas State University, 138 Cardwell Hall, Manhattan, KS 66506}
\email{bergnerj@member.ams.org}

\date{\today}


\begin{abstract}
Much research has been done on structures equivalent to
topological or simplicial groups.  In this paper, we consider
instead simplicial monoids.  In particular, we show that the usual
model category structure on the category of simplicial monoids is
Quillen equivalent to an appropriate model category structure on
the category of simplicial spaces with a single point in degree
zero.  In this second model structure, the fibrant objects are
reduced Segal categories. We then generalize the proof to relate
simplicial categories with a fixed object set to Segal categories
with the same fixed set in degree zero.
\end{abstract}


\maketitle

\section{Introduction}

There has been much work done showing the equivalences of
topological (or simplicial) groups, group-like $A_\infty$-spaces,
and loop spaces.  References for this work include \cite{stash},
\cite{th}, and the more recent \cite{bad} and \cite{bcv}. In this
paper, we show that there are some analogous comparisons to be
made when we work with simplicial monoids rather than simplicial
groups.  We show that, from the perspective of homotopy theory,
simplicial monoids are essentially the same as reduced Segal
categories (Definition \ref{SeCat}).

Reduced Segal categories are a special case of the more general
definition of Segal categories.  These objects, whose name is
suggestive of analogous objects used by Segal in his work on loop
spaces \cite{segal}, were defined (although not by that name) by
Dwyer, Kan, and Smith \cite{dks}.  More recently, they have been
of interest as the 1-dimensional case of the more general notion
of Segal $n$-categories \cite{hs}, and also, as shown by the
author in \cite{thesis}, to be a model for the homotopy theory of
homotopy theories.  The idea behind Segal categories is that they
resemble simplicial categories (or categories enriched over
simplicial sets) but composition of morphisms is only associative
up to higher homotopies.  Theorem \ref{cats} (and, more generally,
\cite[8.6]{thesis}) says that the two are essentially the same, as
least from a homotopy-theoretic viewpoint.

In order to describe the problem more fully, we begin with some
notation and terminology.  Let ${\bf \Delta}$ denote the
cosimplicial category, or category whose objects are finite
ordered sets $[n] = \{0, \ldots ,n\}$ for $n \geq 0$ and whose
morphisms are order-preserving maps between them.  Then $\Deltaop$
is the opposite of this category and is called the simplicial
category. Recall that a simplicial set $X$ is a functor $\Deltaop
\rightarrow \Sets$. We will denote the category of simplicial sets
by $\SSets$. (In the course of this paper we will sometimes refer
to simplicial sets as spaces, due to their homotopy-theoretic
similarity with topological spaces \cite[3.6.7]{hovey}.) Induced
from the maps in $\Deltaop$ are the face maps $d_i:X_n \rightarrow
X_{n-1}$ and degeneracy maps $s_i:X_n \rightarrow X_{n-1}$. A few
of the simplicial sets we will use are the $n$-simplex $\Delta
[n]$ for each $n \geq 0$, its boundary $\dot \Delta [n]$, and the
boundary with the $k$th face removed, $V[n,k]$.  More details
about simplicial sets can be found in \cite[I]{gj}.  We denote by
$|X|$ the topological space given by geometric realization of the
simplicial set $X$ \cite[I.2]{gj}.

More generally, a \emph{simplicial object} in a category $\mathcal
C$ is a functor $\Deltaop \rightarrow \mathcal C$.  In particular,
a functor $\Deltaop \rightarrow \SSets$ is a \emph{simplicial
space} or \emph{bisimplicial set} \cite[IV]{gj}.  Given a
simplicial set $X$, we will also use $X$ to denote the constant
simplicial space with the simplicial set $X$ in each degree.  By
$X^t$ (which should be thought of as ``$X$-transposed") we will
denote the simplicial space such that $(X^t)_n$ is the constant
simplicial set $X_n$, or the simplicial set which has the set
$X_n$ at each level.  We will denote the category of all
simplicial spaces by $\SSets^{\Deltaop}$.

Now, consider a simplicial monoid, by which we mean a simplicial
object in the category $\mathcal Mon$ of monoids, or a functor
$\Deltaop \rightarrow \mathcal Mon$.  We will use an alternate
viewpoint in which we use algebraic theories to define simplicial
monoids.  We begin with the definition of an algebraic theory.
Some references for algebraic theories include chapter 3 of
\cite{bor}, the introduction to \cite{bad}, and section 3 of
\cite{multisort}.

\begin{definition}
An \emph{algebraic theory} $\mathcal T$ is a small category with
finite products and objects denoted $T_n$ for $n \geq 0$. For each
$n$, $T_n$ is equipped with an isomorphism $T_n \cong (T_1)^n$.
Note in particular that $T_0$ is the terminal object in $T$.
\end{definition}

Here we will consider one particular theory, the theory of
monoids, which we denote $\Tm$. To describe this theory, we first
consider the full subcategory of the category of monoids generated
by representatives $T_n$ of the isomorphism classes of free
monoids on $n$ generators. We then define the theory of monoids
$\Tm$ to be the opposite of this category. Thus $T_n$, which is
canonically the coproduct of $n$ copies of $T_1$ in $\mathcal
Mon$, becomes the product of $n$ copies of $T_1$ in $\Tm$.  It
follows that there is a projection map $p_{n,i}: T_n \rightarrow
T_1$ for each $1 \leq i \leq n$ in addition to other monoid maps.
In fact, there are such projection maps in any algebraic theory.
We use them to make the following definition.

\begin{definition} \cite[1.1]{bad}
Given an algebraic theory $\mathcal T$, a \emph{strict simplicial}
$\mathcal T$-\emph{algebra} (or simply $\mathcal
T$-\emph{algebra}) $A$ is a product-preserving functor $A:\mathcal
T \rightarrow \mathcal {SS}ets$.  Here, ``product-preserving"
means that for each $n \geq 0$ the canonical map
\[ A(T_n) \rightarrow A(T_1)^n, \]
induced by the $n$ projection maps $p_{n,i}:T_n \rightarrow T_1$,
is an isomorphism of simplicial sets.  In particular, $A(T_0)$ is
the one-point simplicial set $\Delta [0]$.
\end{definition}

In general, a $\mathcal T$-algebra $A$ defines a strict algebraic
structure on the space $A(T_1)$ corresponding to the theory
$\mathcal T$ \cite[\S 1]{bad}. So, a $\Tm$-algebra $A$ defines a
monoid structure on the space $A(T_1)$.  It turns out that the
category of simplicial monoids is equivalent to the category of
$\Tm$-algebras \cite[II.4]{quillen}.

We can also consider the case where the products are not preserved
strictly, but only up to homotopy.

\begin{definition} \cite[1.2]{bad}
Given an algebraic theory $\mathcal T$, a \emph{homotopy}
$\mathcal T$-\emph{algebra} is a functor $X:\mathcal T \rightarrow
\SSets$ which preserves products up to homotopy.  The functor $X$
preserves products up to homotopy if, for each $n \geq 0$ the
canonical map
\[ X(T_n) \rightarrow X(T_1)^n \]
induced by the projection maps $p_{n,i}:T_n \rightarrow T_1$ is a
weak equivalence of simplicial sets. In particular, we assume that
$X(T_0)$ is weakly equivalent to $\Delta [0]$.
\end{definition}

Now that we have established our definitions for simplicial
monoids, we turn to reduced Segal categories.  These will be
simplicial spaces satisfying some conditions, so, like the
simplicial monoids, they will be given by a diagram of simplicial
sets. We begin with the definition of a Segal precategory.

\begin{definition}
A \emph{Segal precategory} is a simplicial space $X$ such that
$X_0$ is a discrete simplicial set.  If $X_0$ consists of a single
point, then we will call $X$ a \emph{reduced Segal precategory}.
\end{definition}

Now note that for any simplicial space $X$ there is a \emph{Segal
map}
\[ \varphi_k: X_k \rightarrow \underbrace{X_1 \times_{X_0} \ldots
\times_{X_0} X_1}_k \] for each $k \geq 2$, which we define as
follows.  Let $\alpha^i:[1] \rightarrow [k]$ be the map in ${\bf
\Delta}$ such that $\alpha^i(0)=i$ and $\alpha^i(1)=i+1$, defined
for each $0 \leq i \leq k-1$. We can then define the dual map
$\alpha_i:[k]\rightarrow [1]$ in $\Deltaop$.  For $k \geq 2$, the
Segal map is defined to be the map
\[ \varphi_k: X_k \rightarrow \underbrace{X_1 \times_{X_0} \cdots \times_{X_0} X_1}_k \]
induced by the maps
\[ X(\alpha_i):X_k \rightarrow X_1. \]

\begin{definition} \cite[\S 2]{hs} \label{SeCat}
A \emph{Segal category} $X$ is a Segal precategory $X: \Deltaop
\rightarrow \SSets$ such that for each $k \geq 2$ the Segal map
\[ \varphi_k: X_k \rightarrow \underbrace{X_1 \times_{X_0} \cdots \times_{X_0} X_1}_k \]
is a weak equivalence of simplicial sets.  If $X_0= \Delta [0]$,
then $X$ is a \emph{reduced Segal category}.  Note that in this
case we have that $X_k \simeq (X_1)^k$.
\end{definition}

Note that, unlike in the case of homotopy $\mathcal T$-algebras,
we actually require $F([0])$ to be $\Delta [0]$, so that we
actually get a reduced Segal category and not just a simplicial
space equivalent to one.

The precise relationship that we would like to prove between the
simplicial monoids (or $\Tm$-algebras) and reduced Segal
categories is a Quillen equivalence of model categories.  (We will
give a brief overview of model categories and Quillen equivalences
in the next section.)

In this paper, we will prove the existence of a model category
structure $\Algtm$ on the category of $\Tm$-algebras. There is
another model category structure $\LSSp_{*,f}$ on the category of
reduced Segal precategories such that the fibrant objects are
reduced Segal categories.  (We cannot obtain a model category
structure on the category of reduced Segal categories because this
category is not closed under limits and colimits.) Our main result
is the following theorem:

\begin{theorem} \label{main}
The model category structure $\Algtm$ is Quillen equivalent to the
model category structure $\LSSp_{*,f}$.
\end{theorem}

We prove this theorem in section 4 after reviewing model
categories in section 2 and setting up our particular model
category structures in section 3. In section 5, we generalize this
result to an analogous statement for simplicial categories with a
fixed object set $\mathcal O$ and Segal categories with the same
fixed set $\mathcal O$ in degree zero. In \cite{thesis}, we prove
the more general result that there is a Quillen equivalence
between a model category structure on the category of all small
simplicial categories (where the weak equivalences are a
simplicial version of equivalences of categories \cite[\S
1]{simpcat}) and a model category structure on the category of all
Segal precategories (with analogous weak equivalences).

The techniques of this paper carry over to a similar result for
simplicial groups by modifying the category $\Deltaop$.  Let
$I[n]$ denote the category with $n+1$ objects and exactly one
isomorphism between any two objects \cite[\S 6]{rezk}. Similarly
to the way we form ${\bf \Delta}$ from the categories $[n]$, we
define the category ${\bf I\Delta}$ whose objects are the $I[n]$
for all $n \geq 0$ and whose morphisms are generated by the
order-preserving maps and a ``flip" map at each level.

Now, a functor $\IDeltaop \rightarrow \SSets$ is a \emph{reduced
Segal groupoid} if the appropriate analogue to the Segal map is a
weak equivalence for each $k \geq 2$. We can then consider the
theory of groups $\mathcal T_G$ and thus strict $\mathcal
T_G$-algebras, which are essentially simplicial groups. The
arguments of this paper should give us the result that the model
category structure on $\mathcal T_G$-algebras is Quillen
equivalent to an appropriate model structure on ``reduced Segal
pregroupoids."  The explicit argument will be included in a later
paper \cite{inv}.

The idea here, for both monoids and groups, is that when we
consider diagrams of spaces given by a particular algebraic
theory, from a homotopy-theoretic perspective we can actually
consider diagrams of spaces given by some simpler, pared-down
diagram in the sense that $\Deltaop$ can be considered a
simplification of $\mathcal T_M$. Similar work has been done by
Segal \cite{segal} for infinite loop spaces and abelian monoids,
and by Bousfield \cite{bous} for groups and $n$-fold loop spaces.
Berger \cite{berger} has also proved related results for iterated
loop spaces. We hope to find more such examples in the future.  In
particular, it seems that Segal's construction can be modified to
obtain a model for abelian groups \cite{inv}.

\begin{thank}
I am grateful to Bill Dwyer for many helpful conversations about
this paper, Clemens Berger for comments on an earlier draft, and
the referee for suggestions for its improvement.  Early work on
this paper was partially supported by a Clare Boothe Luce
Foundation Graduate Fellowship.
\end{thank}

\section{A Review of Model Categories}

Recall that a \emph{model category structure} on a category
$\mathcal C$ is a choice of three distinguished classes of
morphisms: fibrations, cofibrations, and weak equivalences.  A
(co)fibration which is also a weak equivalence will be called an
\emph{acyclic (co)fibration}.  With this choice of three classes
of morphisms, $\mathcal C$ is required to satisfy five axioms
MC1-MC5 \cite[3.3]{ds}.
%
%
%
%
%
%

An object $X$ in $\mathcal C$ is \emph{fibrant} if the unique map
$X \rightarrow \ast$ from $X$ to the terminal object is a
fibration.  Dually, $X$ is \emph{cofibrant} if the unique map
$\phi \rightarrow X$ from the initial object to $X$ is a
cofibration.  The factorization axiom (MC5) guarantees that each
object $X$ has a weakly equivalent fibrant replacement $\widehat
X$ and a weakly equivalent cofibrant replacement $\widetilde X$.
These replacements are not necessarily unique, but they can be
chosen to be functorial in the cases we will use
\cite[1.1.3]{hovey}.

The model category structures which we will discuss are all
cofibrantly generated.  In a cofibrantly generated model category
$\mathcal C$, there are two sets of morphisms, one of generating
cofibrations and one of generating acyclic cofibrations, such that
a map is a fibration if and only if it has the right lifting
property with respect to the generating acyclic cofibrations, and
a map is an acyclic fibration if and only if it has the right
lifting property with respect to the generating cofibrations
\cite[11.1.2]{hirsch}.  (Recall that if the dotted arrow lift
exists in the solid arrow diagram
\[ \xymatrix{A \ar[r] \ar[d]^i & X \ar[d]^p \\
B \ar[r] \ar@{-->}[ur] & Y} \] then we say that $i$ has the
\emph{left lifting property} with respect to $p$ and that $p$ has
the \emph{right lifting property} with respect to $i$.)  To state
the conditions we will use to determine when we have a cofibrantly
generated model category structure, we make the following
definition.

\begin{definition} \cite[10.5.2]{hirsch}
Let $\mathcal C$ be a category and $I$ a set of maps in $\mathcal
C$.  Then an $I$-\emph{injective} is a map which has the right
lifting property with respect to every map in $I$.  An
$I$-\emph{cofibration} is a map with the left lifting property
with respect to every $I$-injective.
\end{definition}

Using this definition, we will use the following theorem to
determine when we have a cofibrantly generated model category
structure.

\begin{theorem} \cite[11.3.1]{hirsch} \label{CofGen}
Let $\mathcal M$ be a category which has all small limits and
colimits. Suppose that $\mathcal M$ has a class of weak
equivalences which satisfies the two-out-of-three property (model
category axiom MC2) and which is closed under retracts. Let $I$
and $J$ be sets of maps in $\mathcal M$ which satisfy the
following conditions:
\begin{enumerate}
\item Both $I$ and $J$ permit the small object argument
\cite[10.5.15]{hirsch}.

\item Every $J$-cofibration is an $I$-cofibration and a weak
equivalence.

\item Every $I$-injective is a $J$-injective and a weak
equivalence.

\item One of the following conditions holds:
\begin{enumerate}
\item A map that is an $I$-cofibration and a weak equivalence is a
$J$-cofibration, or

\item A map that is both a $J$-injective and a weak equivalence is
an $I$-injective.
\end{enumerate}
\end{enumerate}
Then there is a cofibrantly generated model category structure on
$\mathcal M$ in which $I$ is a set of generating cofibrations and
$J$ is a set of generating acyclic cofibrations.
\end{theorem}

We will refer to the standard model category structure on the
category $\SSets$ of simplicial sets.  In this case, a weak
equivalence is a map of simplicial sets $f:X \rightarrow Y$ such
that the induced map $|f|:|X| \rightarrow |Y|$ is a weak homotopy
equivalence of topological spaces.  The cofibrations are
monomorphisms, and the fibrations are the maps with the right
lifting property with respect to the acyclic cofibrations
\cite[I.11.3]{gj}. This model category structure is cofibrantly
generated; the generating cofibrations are the maps $\dot \Delta
[n] \rightarrow \Delta [n]$ for $n \geq 0$, and the generating
acyclic cofibrations are the maps $V[n,k] \rightarrow \Delta [n]$
for $n \geq 1$ and $0 \leq k \leq n$.

Each model structure in this paper will have the additional
structure of a simplicial model category.  A \emph{simplicial
category} is a category with a simplicial set of morphisms between
any two objects.  A \emph{simplicial model category}, then, is a
model category $\mathcal M$ which is also a simplicial category,
satisfying some axioms \cite[9.1.6]{hirsch}. (Notice that this
terminology is somewhat confusing in that simplicial objects in
the category of categories are not necessarily simplicial
categories, nor are simplicial objects in the category of model
categories actually simplicial model categories, at least not
without imposing further conditions.) The important part of this
structure is that it enables us to talk about the \emph{function
complex}, or simplicial set $\Map(X,Y)$ for any pair of objects
$X$ and $Y$ in $\mathcal M$.  For example, the model category
$\SSets$ is a simplicial model category.  For simplicial sets $X$
and $Y$, the simplicial set $\Map(X,Y)$ is given by
\[ \Map(X,Y)_n = \Hom (X \times \Delta [n], Y). \]
However, a function complex is not necessarily homotopy invariant,
so we have the following definition.

\begin{definition} \cite[17.3.1]{hirsch}
If $\mathcal M$ is a simplicial model category, then the
\emph{homotopy function complex} $\Map^h(X,Y)$ is given by
$\Map(\widetilde X, \widehat Y)$, where $\widetilde X$ is a
cofibrant replacement of $X$ in $\mathcal M$ and $\widehat Y$ is a
fibrant replacement for $Y$.
\end{definition}

More generally, a homotopy function complex $\Map^h(X,Y)$ in a
(not necessarily simplicial) model category $\mathcal M$, or in a
category with weak equivalences, is the simplicial set $\Map(X,Y)$
given by the simplicial set of morphisms between $X$ and $Y$ in
the hammock localization $L^H \mathcal M$ of $\mathcal M$
\cite{dk1}.

Many of the model category structures that we will work with will
be obtained by localizing a given model category structure with
respect to a set of maps.  To make sense of this notion, we make
the following definitions.

\begin{definition} \label{local}
Suppose that $S = \{f:A \rightarrow B\}$ is a set of maps in a
model category $\mathcal M$.  An $S$-\emph{local} object $X$ is a
fibrant object of $\mathcal M$ such that for any $f:A \rightarrow
B$ in $S$, the induced map on homotopy function complexes
\[ f^*:\Map^h(B,W) \rightarrow \Map^h(A,W) \]
is a weak equivalence of simplicial sets.  A map $g:X \rightarrow
Y$ in $\mathcal M$ is then an $S$-\emph{local equivalence} if for
every local object $W$, the induced map on homotopy function
complexes
\[ g^*: \Map^h(Y,W) \rightarrow \Map^h(X,W) \]
is a weak equivalence of simplicial sets.
\end{definition}

We now define the localization of a model category $\mathcal M$.
The following theorem holds for all model categories $\mathcal M$
which are left proper and cellular. We will not define these
properties here but refer the reader to \cite[13.1.1,
12.1.1]{hirsch} for details.  All the model categories we will
work with can be shown to satisfy these properties.

\begin{theorem}\cite[4.1.1]{hirsch} \label{Loc}
Let $\mathcal M$ be a left proper cellular model category.  There
is a model category structure $\mathcal L_S \mathcal M$ on the
underlying category of $\mathcal M$ such that:
\begin{enumerate}
\item The weak equivalences are the $S$-local equivalences.

\item The cofibrations are precisely the cofibrations of $\mathcal
M$.

\item The fibrations are the maps which have the right lifting
property with respect to the maps which are both cofibrations and
$S$-local equivalences.

\item The fibrant objects are the $S$-local objects.

\item If $\mathcal M$ is a simplicial model category, then its
simplicial structure induces a simplicial structure on $\mathcal
L_S \mathcal M$.
\end{enumerate}
\end{theorem}

In particular, given an object $X$ of $\mathcal M$, we can talk
about its functorial fibrant replacement $LX$ in $\mathcal L_S
\mathcal M$. The object $LX$ is an $S$-local object which is
fibrant in $\mathcal M$, and we will call it the
\emph{localization} of $X$ in $\mathcal L_S \mathcal M$.

We now state the definition of a Quillen pair of model category
structures.  Recall that for categories $\mathcal C$ and $\mathcal
D$ a pair of functors
\[ \xymatrix@1{F: \mathcal C \ar@<.5ex>[r] & \mathcal D:R
\ar@<.5ex>[l]} \] is an \emph{adjoint pair} if for each object $X$
of $\mathcal C$ and object $Y$ of $\mathcal D$ there is a natural
isomorphism $\psi :\Hom_\mathcal D(FX,Y) \rightarrow \Hom_\mathcal
C(X,RY)$ \cite[IV.1]{macl}.  The adjoint pair is sometimes written
as the triple $(F, R, \psi)$.

\begin{definition} \cite[1.3.1]{hovey}
If $\mathcal C$ and $\mathcal D$ are model categories, then an
adjoint pair $(F,R, \psi)$ between them is a \emph{Quillen pair}
if one of the following equivalent statements holds:
\begin{enumerate}
\item $F$ preserves cofibrations and acyclic cofibrations.

\item $R$ preserves fibrations and acyclic fibrations.
\end{enumerate}
\end{definition}

We will use the following theorem to show that we have a Quillen
pair of localized model category structures.

\begin{theorem} \cite[3.3.20]{hirsch} \label{LocPair}
Let $\mathcal C$ and $\mathcal D$ be left proper, cellular model
categories and let $(F,R, \psi)$ be a Quillen pair between them.
Let $S$ be a set of maps in $\mathcal C$ and $L_S \mathcal C$ the
localization of $\mathcal C$ with respect to $S$.  Then if ${\bf
L}FS$ is the set in $\mathcal D$ obtained by applying the left
derived functor of $F$ to the set $S$ \cite[8.5.11]{hirsch}, then
$(F,R, \psi)$ is also a Quillen pair between the model categories
$L_S\mathcal C$ and $L_{{\bf L}FS}D$.
\end{theorem}

We now have the following definition of Quillen equivalence, which
is the standard notion of ``equivalence" of model category
structures.

\begin{definition} \cite[1.3.12]{hovey}
A Quillen pair is a \emph{Quillen equivalence} if for all
cofibrant $X$ in $\mathcal C$ and fibrant $Y$ in $\mathcal D$, a
map $f:FX \rightarrow Y$ is a weak equivalence in $\mathcal D$ if
and only if the map $\varphi f:X \rightarrow RY$ is a weak
equivalence in $\mathcal C$.
\end{definition}

We will use the following proposition to prove that our Quillen
pairs are Quillen equivalences.  Recall that a functor $F:\mathcal
C \rightarrow \mathcal D$ \emph{reflects} a property if, for any
morphism $f$ of $\mathcal C$, whenever $Ff$ has the property, then
so does $f$.

\begin{prop}\cite[1.3.16]{hovey}
Suppose that $(F,R,\psi)$ is a Quillen pair from $\mathcal C$ to
$\mathcal D$. Then the following statements are equivalent:
\begin{enumerate}
\item $(F,R, \psi)$ is a Quillen equivalence.

\item $F$ reflects weak equivalences between cofibrant objects,
and for every fibrant $Y$ in $\mathcal D$ the map $F(RY)^c
\rightarrow Y$ is a weak equivalence.

\item $R$ reflects weak equivalences between fibrant objects, and
for every cofibrant $X$ in $\mathcal C$ the map $X \rightarrow
R(FX)^f$ is a weak equivalence.
\end{enumerate}
\end{prop}

\section{Model Category Structures}

In this section, we set up the model category structures that we
will need to prove Theorem \ref{main}.

Because we will need to consider a model structure on the category
of simplicial monoids, we now describe a model structure on the
category of $\mathcal T$-algebras for any algebraic theory
$\mathcal T$.

\begin{prop} \cite[II.4]{quillen}, \cite[3.1]{sch}
Let $\mathcal T$ be an algebraic theory and $\Algt$ the category
of $\mathcal T$-algebras.  Then there is a cofibrantly generated
model category structure on $\Algt$ in which the weak equivalences
and fibrations are levelwise weak equivalences of simplicial sets
and the cofibrations are the maps with the left lifting property
with respect to the maps which are fibrations and weak
equivalences.
\end{prop}


We would also like to have a model category structure for homotopy
$\mathcal T$-algebras. However, we will not have a model category
structure on the category of homotopy $\mathcal T$-algebras
itself; the category of homotopy $\mathcal T$-algebras does not
have all coproducts. We will have a model category structure on
the category of all $\mathcal T$-diagrams of simplicial sets in
which homotopy $\mathcal T$-algebras are the fibrant objects.  To
obtain this structure, we begin by considering the category of all
functors $\mathcal T \rightarrow \SSets$, which we denote by
$\SSetst$.

\begin{theorem} \cite[IX 1.4]{gj}
There is a model category structure on $\SSetst$ in which the weak
equivalences are the levelwise weak equivalences, the fibrations
are the levelwise fibrations, and the cofibrations are the maps
which have the left lifting property with respect to the maps
which are both fibrations and weak equivalences.
\end{theorem}

The desired model structure can be obtained by localizing the
model structure on $\SSetst$ with respect to a set of maps. We
will summarize this localization here; a complete description is
given by Badzioch \cite[\S 5]{bad}.

Given an algebraic theory $\mathcal T$, consider the functor
$\Hom_\mathcal T (T_k,-)$.  We then have maps
\[ p_k: \coprod_{i=1}^k \Hom_\mathcal T(T_1,-) \rightarrow \Hom_\mathcal T(T_k,-) \]
induced from the projection maps in $\mathcal T$.  We then
localize the model category structure on $\SSetst$ with respect to
the set $S= \{p_k | k \geq 0 \}$. We denote the resulting model
category structure $\LSSetst$.

\begin{prop} \cite[5.5]{bad}.
The fibrant objects in $\LSSetst$ are the homotopy $\mathcal
T$-algebras which are fibrant in $\SSetst$.
\end{prop}

We now have the following result by Badzioch.

\begin{theorem}  \cite[6.4]{bad} \label{rigid}
Given an algebraic theory $\mathcal T$, there is a Quillen
equivalence of model categories between $\Algt$ and $\LSSetst$.
\end{theorem}

By applying this Quillen equivalence to the theory of monoids
$\Tm$, we have reduced the proof of our main theorem to finding a
Quillen equivalence between $\LSSetstm$ and an appropriate model
structure for reduced Segal categories.  As with the homotopy
$\mathcal T$-algebras, we will need to consider the category of
reduced Segal precategories and find a model structure in which
the fibrant objects are Segal categories.

We consider the general case of the category whose objects are
Segal precategories with a fixed set $\mathcal O$ in degree zero
and whose morphisms are object-preserving on that set.  We will
prove the existence of a model structure on this category in which
the weak equivalences are levelwise weak equivalences of
simplicial sets. We give the proof here for any set $\mathcal O$,
since we will need it for the more general construction in the
last section of the paper. Thus, let $\SSpo$ denote the category
whose objects are simplicial spaces with a fixed set $\mathcal O$
in degree zero and whose morphisms are maps of simplicial spaces
which are the identity map on $\mathcal O$.

This model category structure on $\SSpo$ will be analogous to the
projective model category structure on simplicial spaces in which
the fibrations and weak equivalences are levelwise fibrations and
weak equivalences of simplicial sets \cite[IX.1.4]{gj}. However,
this structure must be modified so that the objects defining the
generating cofibrations and generating acyclic cofibrations will
be in the category $\SSpo$. We begin by understanding limits and
colimits in $\SSpo$.

\begin{lemma} \label{lim}
$\SSpo$ has all small limits.
\end{lemma}

\begin{proof}
Suppose that for each object $\alpha$ of an index category
$\mathcal D$, $X_\alpha$ is a simplicial space with the set
$\mathcal O$ in degree zero. In the category of all simplicial
spaces, we have the limit $\lim_\mathcal D X_\alpha$. However, if
$\mathcal D$ is not connected, then this limit will not be in the
category $\SSpo$ since $\lim_\mathcal D (X_\alpha)_0 = \mathcal
O^{\pi_0 \mathcal D}$. However, if $\text{diag}:\mathcal O
\rightarrow \lim_\mathcal D (X_\alpha)_0$ is the diagonal map, we
can define the limit in $\SSpo$, denoted $\lim_\mathcal D^\mathcal
O X_\alpha$, as the pullback in the diagram
\[ \xymatrix{\lim_\mathcal D^\mathcal O X_\alpha \ar[r] \ar[d] &
\lim_\mathcal D X_\alpha \ar[d] \\
\mathcal O \ar[r]^-{\text{diag}} & \lim_\mathcal D(X_\alpha)_0} \]
This new object now satisfies the universal property of limits
when we require the maps involved to be the identity on degree
zero and hence in our category $\SSpo$.
\end{proof}

\begin{lemma} \label{colim}
$\SSpo$ has all small colimits.
\end{lemma}

\begin{proof}
As with the limits, begin by considering the ordinary colimit in
simplicial spaces.  Again, let $X_\alpha$, indexed by the objects
$\alpha$ of $\mathcal D$, be objects in $\SSpo$.  Note that again
we have a problem in degree zero if our index category $\mathcal
D$ has more than one component, since in this case
$\colim_\mathcal D (X_\alpha)_0 = \coprod_{\pi_0 \mathcal D}
\mathcal O$. Then if we consider the fold map
\[ \text{fold}:\colim_\mathcal D(X_\alpha)_0
\rightarrow \mathcal O, \] we can define the colimit in $\SSpo$,
denoted $\colim_D^\mathcal O X_\alpha$ as the pushout in the
diagram
\[ \xymatrix{\colim_\mathcal D (X_\alpha)_0 \ar[r]^-{\text{fold}}
\ar[d] & \mathcal O \ar[d] \\
\colim_\mathcal D X_\alpha \ar[r] & \colim_\mathcal D^\mathcal O
X_\alpha} \] where the left-hand vertical map is the inclusion
map.  Similarly to the case for limits, this new simplicial space
satisfies the universal property for colimits.
\end{proof}

\begin{prop} \label{sspof}
There is a model category structure on $\SSpo$, which we denote
$\SSpof$, in which the weak equivalences are levelwise weak
equivalences of simplicial sets, the fibrations are the levelwise
fibrations of simplicial sets, and the cofibrations are the maps
with the left lifting property with respect to the maps which are
fibrations and weak equivalences.
\end{prop}

The first step towards defining this model category structure is
finding candidates for the generating cofibrations and generating
acyclic cofibrations.  Recall that for the projective model
category structure on the category $\SSets^{\Deltaop}$ of
simplicial spaces, in which the weak equivalences and fibrations
are levelwise, the generating acyclic cofibrations are of the form
\[ V[m,k] \times \Delta [n]^t \rightarrow \Delta [m] \times \Delta [n]^t \]
\cite[IV.3.1]{gj}. (Recall that by $\Delta [n]^t$ we denote the
simplicial space which is the constant simplicial set $\Delta
[n]_k$ in degree $k$.)

The first problem in using these maps for $\SSpof$ is that $\Delta
[n]^t$ is not going to be in $\SSpo$ for all values of $n$.
Instead, we need to define a separate $n$-simplex for any
$n$-tuple $x_0, \ldots ,x_n$ of objects in $\mathcal O$, denoted
$\Delta [n]^t_{x_0, \ldots ,x_n}$, so that the objects are
preserved. Setting $\xu =(x_1, \ldots ,x_n)$, we write this
simplex as $\Delta [n]^t_\xu$.  (The values of the $x_i$ can
repeat in a particular $\xu$.)  Note that this object $\Delta
[n]^t_\xu$ also needs to have all elements of $\mathcal O$ as
0-simplices, so we add any of these elements that have not already
been included in the $x_i$'s, plus their degeneracies in higher
degrees.

Now we consider the modified generating acyclic cofibrations
\[ V[m,k] \times \Delta [n]^t_\xu \rightarrow \Delta
[m] \times \Delta [n]^t_\xu. \]  However, these objects are still
not in $\SSpof$ because the simplicial set in degree zero for each
is not discrete. This problem can be fixed by collapsing the
copies of $\Delta [m]$ and $V[m,k]$ in degree zero to their
respective $x_i$'s. Note that the degeneracies then get collapsed
as well. More explicitly, we define the object $(R_{m,n,k})_\xu$
to be the pushout of the diagram
\[ \xymatrix{V[m,k] \times (\Delta [n]^t_\xu)_0 \ar[r] \ar[d] & V[m,k]
\times \Delta [n]^t_\xu \ar[d] \\
(\Delta [n]^t_\xu)_0 \ar[r] & (R_{m,n,k})_\xu.} \] Similarly, we
define the object $(Q_{m,n})_\xu$ to be the pushout of the diagram
\[ \xymatrix{\Delta[m] \times (\Delta [n]^t_\xu)_0 \ar[r] \ar[d] &
\Delta [m] \times \Delta [n]^t_\xu \ar[d] \\
(\Delta [n]^t_\xu)_0 \ar[r] & (Q_{m,n})_\xu.} \]

Now we are able to define the set of maps
\[ J_f= \{ (R_{m,n,k})_\xu \rightarrow (Q_{m,n})_\xu \} \] where $m,n \geq 1$, $0 \leq k
\leq m$, and $\xu =(x_0, \ldots ,x_n) \in \mathcal O^{n+1}$. This
set $J_f$ will be a set of generating acyclic cofibrations for
$\SSpof$.

Similarly, we can define the set
\[ I_f= \{(P_{m,n})_\xu \rightarrow (Q_{m,n})_\xu \} \] for all $m,n \geq 0$ and $\xu \in \mathcal O^{n+1}$, where
$(P_{m,n})_\xu$ is the pushout of the diagram
\[ \xymatrix{\dot \Delta [m] \times (\Delta [n]^t_\xu)_0 \ar[d] \ar[r] &
\dot \Delta [m] \times \Delta [n]^t_\xu \ar[d] \\
(\Delta [n]^t_\xu)_0 \ar[r] & (P_{m,n})_\xu.} \]  We will show
that these maps are a set of generating cofibrations for $\SSpof$.

\begin{proof} [Proof of Theorem \ref{sspof}]
Lemmas \ref{lim} and \ref{colim} show that our category has small
limits and colimits. The two-out-of-three property and the retract
axiom for weak equivalences follow as usual; see for example
\cite[8.10]{ds}. It now suffices to check the conditions of
Theorem \ref{CofGen}. To prove statement (1), notice that the maps
in the sets $I_f$ and $J_f$ are modified versions of the
generating cofibrations in the projective model category structure
on simplicial spaces, which permit the small object argument
\cite[10.5.15]{hirsch}.  Hence the ones in $I_f$ and $J_f$ do
also.

Notice that the $I_f$-injectives are precisely the levelwise
acyclic fibrations, and that the $J_f$-injectives are precisely
the levelwise fibrations.  Thus, we have satisfied conditions (3)
and (4)(ii).

Now notice that the $I_f$-cofibrations are precisely the
cofibrations, by our definition of cofibration. Furthermore, the
$J_f$-cofibrations are the maps with the left lifting property
with respect to the fibrations.
It follows that a $J_f$-cofibration is an $I_f$-cofibration. Using
the definition of the maps in $J_f$ and simplicial set arguments,
we can see that a $J_f$-cofibration is also a weak equivalence.
\end{proof}

We now want to localize this model category so that the fibrant
objects are Segal categories with the set $\mathcal O$ in degree
zero. If we can find an appropriate map to localize with respect
to, then the desired model category structure will follow from
Theorem \ref{Loc}.

We first consider the map $\varphi$ used by Rezk \cite[\S 4]{rezk}
to localize simplicial spaces to obtain more general Segal spaces,
then modify it so that the objects are in $\SSpof$.


Rezk defines a map $\alpha^i:[1] \rightarrow [k]$ in ${\bf
\Delta}$ such that $0 \mapsto i$ and $1 \mapsto i+1$ for each $0
\leq i \leq k-1$. Then for each $k$ he defines the object
\[ G(k)^t= \bigcup_{i=0}^{k-1} \alpha^i \Delta [1]^t \] and the inclusion map $\varphi^k:
G(n)^t \rightarrow \Delta [k]^t$.  His localization is with
respect to the coproduct of inclusion maps
\[ \varphi = \coprod_{k \geq 0} (G(k)^t \rightarrow \Delta [k]^t). \]

However, in our case, the objects $G(k)^t$ and $\Delta [k]^t$ are
not in the category $\SSpo$.  As before, we can replace $\Delta
[k]^t$ with the objects $\Delta [k]^t_\xu$, where $\xu =(x_0,
\ldots ,x_k)$. Then, define
\[ G(k)^t_\xu = \bigcup_{i=0}^{k-1} \alpha^i \Delta
[1]^t_{x_i, x_{i+1}}. \] Now, we need to take coproducts not only
over all values of $k$, but also over all $k$-tuples of vertices.

We first define for each $k \geq 0$ the map
\[ \varphi^k = \coprod_{\xu \in \mathcal O^{k+1}}(G(k)^t_\xu
\rightarrow \Delta [k]^t_\xu). \]  Then the map $\varphi$ looks
like
\[ \varphi = \coprod_{k \geq 0}(\varphi^k: \coprod_{\xu \in
\mathcal O^{k+1}} (G(k)^t_\xu \rightarrow \Delta [k]^t_\xu)).
\] When the set $\mathcal O$ is not clear from the context, we
will write $\varphi_\mathcal O$ to specify that we are in $\SSpo$.

For any simplicial space $X$, there is a weak equivalence of
simplicial sets
\[ \Map^h(\coprod_{\xu \in \mathcal O^{k+1}}G(k)^t_\xu ,X) \rightarrow \underbrace{X_1 \times^h_{X_0} \cdots
\times^h_{X_0} X_1}_k \] where the right-hand side is the homotopy
limit of the diagram
\[ \xymatrix{X_1 \ar[r]^{d_0} & X_0 & X_1 \ar[l]_{d_1}
\ar[r]^{d_0} & \ldots \ar[r]^{d_0} & X_0 & X_1 \ar[l]_{d_1}} \]
with $k$ copies of $X_1$.  However, in this case, since $X_0$ is
discrete we can take the limit
\[ \underbrace{X_1 \times_{X_0} \cdots \times_{X_0} X_1}_k \] on
the right hand side.

Then for any simplicial space $X$ there is a map
\[ \varphi_k = \Map^h(\varphi^k, X): \Map^h(\coprod_{\xu} \Delta
[k]^t_\xu,X) \rightarrow \Map^h(\coprod_{\xu} G(k)^t_\xu ,X). \]
More simply written, this map is
\[ \varphi_k:X_k \rightarrow \underbrace{X_1 \times_{X_0} \cdots
\times_{X_0} X_1}_k \] and id precisely the Segal map used in the
definition of a Segal category.


\begin{prop}
Localizing the model category structure on $\SSpof$ with respect
to the map $\varphi_\mathcal O$ results in a model category
structure $\LSSpof$ on simplicial spaces with a fixed set
$\mathcal O$ in degree zero in which the weak equivalences are the
$\varphi_\mathcal O$-local equivalences, the cofibrations are
those of $\SSpof$, and the fibrations are the maps with the right
lifting property with respect to the cofibrations which are
$\varphi_\mathcal O$-local equivalences.
\end{prop}

\begin{proof}
The proof follows from Theorem \ref{Loc}.
\end{proof}

In \cite{fibrant}, we prove that the fibrant objects in $\LSSpof$
are precisely the Segal categories with $\mathcal O$ in degree
zero which are fibrant in the projective model category structure.

For making some of our calculations, we will find it convenient to
work in a model category structure $\SSpoc$ in which the weak
equivalences are again given by levelwise weak equivalences of
simplicial sets, but in which the cofibrations, rather than the
fibrations, are levelwise.

\begin{theorem} \label{sspoc}
There is a model category structure $\SSpoc$ on the category of
Segal precategories with a fixed set $\mathcal O$ in degree zero
in which the weak equivalences and cofibrations are levelwise, and
in which the fibrations are the maps with the right lifting
property with respect to the acyclic cofibrations.
\end{theorem}

To define sets $I_c$ and $J_c$ which will be our candidates for
generating cofibrations and generating acyclic cofibrations,
respectively, we first recall the generating cofibrations and
acyclic cofibrations in the Reedy (or injective) model category
structure on simplicial spaces, in which the weak equivalences and
cofibrations are levelwise. The generating cofibrations are the
maps
\[ \dot \Delta [m] \times \Delta [n]^t \cup \Delta [m] \times \dot
\Delta [n]^t \rightarrow \Delta [m] \times \Delta [n]^t \] for all
$m,n \geq 0$, and similarly the generating acyclic cofibrations
are the maps
\[ V[m,k] \times \Delta [n]^t \cup \Delta [m] \times \dot \Delta
[n]^t \rightarrow \Delta [m] \times \Delta [n]^t \] for all $n
\geq 0$, $m \geq 1$, and $0 \leq k \leq m$ \cite[2.4]{rezk}.


To modify these maps, we begin by considering the category
$\Secat$ of all Segal precategories and the inclusion functor
$\Secat \rightarrow \SSets^{\Deltaop}$. This functor has a left
adjoint which we call the reduction functor. Given a simplicial
space $X$, we denote its reduction by $(X)_r$. Reducing $X$
essentially amounts to collapsing the space $X_0$ to its set of
components and making the appropriate changes to degeneracies in
higher degrees.  So, we start by reducing the objects defining the
Reedy generating cofibrations and generating acyclic cofibrations
to obtain maps of the form
\[ (\dot \Delta [m] \times \Delta [n]^t \cup \Delta [m] \times \dot
\Delta [n]^t)_r \rightarrow (\Delta [m] \times \Delta [n])_r \]
and
\[ (V[m,k] \times \Delta [n]^t \cup \Delta [m] \times \dot \Delta
[n]^t)_r \rightarrow (\Delta [m] \times \Delta [n]^t)_r \] Then,
in order to have our maps in $\SSpo$, we define a separate such
map for each choice of vertices $\xu$ in degree zero and adding in
the remaining points of $\mathcal O$ if necessary.  As above, we
use $\Delta [n]^t_\xu$ to denote the object $\Delta [n]^t$ with
the $(n+1)$-tuple $\xu$ of vertices. We then define sets
\[ I_c = \{(\dot \Delta [m] \times \Delta [n]^t_\xu \cup \Delta
[m] \times \dot \Delta [n]^t_\xu)_r \rightarrow (\Delta [m] \times
\Delta [n]^t_\xu)_r \} \] for all $m \geq 0$ and $n \geq 1$, and
\[ J_c = \{(V[m,k] \times \Delta [n]^t_\xu \cup \Delta [m] \times
\dot \Delta [n]^t_\xu)_r \rightarrow (\Delta [m] \times \Delta
[n]^t_\xu)_r\} \] for all $m \geq 1$, $n \geq 1$, and $0 \leq k
\leq m$, where the notation $(-)_\xu$ indicates the vertices.

Given these maps, we are now able to prove the existence of the
model category structure $\SSpoc$.

\begin{proof} [Proof of Theorem \ref{sspoc}]
The proofs that $\SSpoc$ has finite limits and colimits and
satisfies the two out of three and retract axioms, as well as
condition (1), are the same as for $\SSpof$.

From the definitions of $I_c$ and $J_c$, it follows that the
$I_c$-injectives are the maps of Segal precategories with
$\mathcal O$ in degree zero which are Reedy fibrations and that
the $J_c$-injectives are the maps of Segal precategories with
$\mathcal O$ in degree zero which are Reedy acyclic fibrations.
Furthermore, it follows from these facts that the
$I_c$-cofibrations are precisely the cofibrations and that the
$J_c$-cofibrations are precisely the acyclic cofibrations.  These
observations imply that the conditions of Theorem \ref{CofGen} are
satisfied.
\end{proof}

We can then localize $\SSpoc$ with respect to the map
$\varphi_\mathcal O$ to obtain a model category structure which we
denote $\LSSpoc$.

At first glance, one might wonder if the weak equivalences in
$\LSSpof$ and $\LSSpoc$ are actually the same, since each of these
model structures is obtained via localization of a different model
structure on Segal precategories. However, since the weak
equivalences before localization are the same in each case, this
localization is independent of the underlying model structure
\cite[\S 7]{thesis}.  Therefore, the weak equivalences are
actually the same in these two localized structures.  (Notice that
here we need to use the more general notion of homotopy function
complex in a category with specified weak equivalences, rather
than in a simplicial model category, since we are working in two
different model categories on the same underlying category.)

We then have the following result.

\begin{prop}
The adjoint pair given by the identity functor induces a Quillen
equivalence of model categories
\[ \xymatrix@1{\LSSpof \ar@<.5ex>[r] & \LSSpoc. \ar@<.5ex>[l]} \]
\end{prop}

\begin{proof}
Since the cofibrations in $\LSSpof$ are monomorphisms, the
identity functor
\[ \LSSpof \rightarrow \LSSpoc \] preserves both cofibrations
and acyclic cofibrations, so this adjoint pair is a Quillen pair.
It remains to show that for any cofibrant $X$ in $\LSSpof$ and
fibrant $Y$ in $\LSSpoc$, the map $FX \rightarrow Y$ is a weak
equivalence if and only of the map $X \rightarrow RY$ is a weak
equivalence.  However, this fact follows since the weak
equivalences are the same in each category.
\end{proof}

Before proceeding to the proof of the main result, we need one
more model structure. Let $\ast$ denote the set with one element.
The objects of $\LSSp_{*,f}$ have a single point in degree zero,
whereas the objects in $\LSSetstm$ have an arbitrary simplicial
set in degree zero.  To simplify matters, we define a model
structure analogous to $\LSSetstm$ but on the category of functors
$\Tm \rightarrow \Sets$ which send $T_0$ to $\Delta [0]$.

\begin{prop}
Consider the category $\SSetstm_*$ of functors $\Tm \rightarrow
\SSets$ such that the image of $T_0$ is $\Delta [0]$.  There is a
model category structure on $\SSetstm_*$ in which the weak
equivalences and fibrations are defined levelwise and the
cofibrations are the maps with the left lifting property with
respect to the acyclic fibrations.
\end{prop}

\begin{proof}
Limits and colimits exist in $\SSetstm_*$ by analogous arguments
to the ones given in Lemmas \ref{lim} and \ref{colim}.  By taking
sets of generating cofibrations and generating acyclic
cofibrations for $\SSetstm$ and modifying them in the same way as
we did for $\SSp_{*,f}$, we can obtain sets of generating
cofibrations and generating acyclic cofibrations for $\SSetstm_*$.
The proof follows just as the proof of Theorem \ref{sspof}.
\end{proof}

Now, to obtain a localized model category $\LSSetstm_*$, we need
to modify the maps
\[ p_k: \coprod_{i=1}^k \Hom_\mathcal T(T_1,-) \rightarrow \Hom_\mathcal
T(T_k,-) \] that we used to obtain $\LSSetstm$ from $\SSetstm$.
Since $\Hom_\mathcal T (T_k, T_0)$ is a single point for all $k$,
the only change we need to make is to take the coproduct
$\coprod_k \Hom_\mathcal T(T_1,-)$ in the category $\SSetstm_*$
(as in Lemma \ref{colim}). We then localize $\SSetstm_*$ with
respect to the set of all such maps to obtain a model structure
$\LSSetstm_*$.

Since a fibrant and cofibrant object $X$ in $\LSSetstm$ has $X_0$
weakly equivalent to $\Delta [0]$, it is not too surprising that
we have the following result:

\begin{prop} \label{tmmcs}
There is a Quillen equivalence of model categories
\[ \xymatrix@1{I:\LSSetstm_* \ar@<.5ex>[r] & \LSSetstm:R. \ar@<.5ex>[l]} \]
\end{prop}

The left adjoint $I: \LSSetstm_* \rightarrow \LSSetstm$ is just
the inclusion functor.  It should be noted that the right adjoint
functor $R$ is not the usual ``collapse" functor, which would be
the left adjoint to $I$. Given an object $Y$ of $\LSSetstm$, $RY$
is given by the pullback
\[ \xymatrix{RY \ar[r] \ar[d] & \Delta[0] \ar[d] \\
Y \ar[r] & Y_0} \] (where $\Delta[0]$ and $Y_0$ denote constant
diagrams).  A detailed description of an analogous functor can be
found in \cite[6.1]{thesis}.

\begin{proof}
The right adjoint $R$ sends the space in degree zero to a point.
Since the cofibrations of $\LSSetstm_*$ are also cofibrations in
$\LSSetstm$ and the weak equivalences are defined in the same way
in each model structure, $I$ preserves cofibrations and acyclic
cofibrations.  Hence, this adjoint pair is a Quillen pair.

By the same argument, $I$ also reflects weak equivalences between
cofibrant objects.  To prove that we have a Quillen equivalence,
it remains to show that for any fibrant object $Y$ in $\LSSetstm$,
the map $I(RY)^c \rightarrow Y$ is a weak equivalence.  Since $I$
is an inclusion, it suffices to show that the map $RY \rightarrow
Y$ is a weak equivalence whenever $Y$ is fibrant.

Since $Y$ is fibrant, and so a homotopy $\Tm$-algebra, we have
that the map $\Delta[0] \rightarrow Y_0$ is a weak equivalence.
Since $\Delta[0]$ and $Y_0$ are fibrant, this map is a levelwise
weak equivalence of simplicial sets \cite[3.2.18]{hirsch}.
Furthermore, the map $Y \rightarrow Y_0$ is a fibration, again
given levelwise. Since the pullback is also defined levelwise, the
fact that $\SSets$ is right proper (i.e., every pullback of a weak
equivalence along a fibration is a weak equivalence), we can
conclude that the map $RY \rightarrow Y$ is a weak equivalence.
\end{proof}

\section{Proof of Theorem \ref{main}}

For the rest of this section, we will set $\mathcal O = *$, the
one-element set. In the previous section, we used Theorem
\ref{rigid} to see that the model category structure on $\Algtm$,
which is equivalent to the category of simplicial monoids, is
Quillen equivalent to $\LSSetstm$. This model category is in turn
Quillen equivalent to $\LSSetstm_*$ by Proposition \ref{tmmcs}.
So, to prove Theorem \ref{main} it suffices to show that there is
a Quillen equivalence between the model categories
$\mathcal{LSS}p_{*,f}$ and $\LSSetstm_*$.

Let $L_1$ be the functorial fibrant replacement functor (or
localization) for $\LSSp_*$ given by taking a colimit of pushouts
along the generating acyclic cofibrations, and let $L_2$ be the
analogous fibrant replacement functor for $\LSSetstm_*$.

We would like to know what the localization $L_1$ does to an
$n$-simplex $\Delta [n]^t_*$.  Note that we will use this notation
in the sense of the previous section: we are specifying that we
are working in the category of reduced Segal precategories, namely
in the case where $\mathcal O=*$. To make the calculations about
our localizations, we will use the model structure
$\LSSp_{\ast,c}$, since in this case all objects are cofibrant and
in particular all monomorphisms are cofibrations.

\begin{lemma} \label{Lhocolim}
Let $L$ be a localization functor for a model category $\mathcal
M$. Given a small diagram of objects $X_\alpha$ of $\mathcal M$,
\[ L(\hocolim X_\alpha) \simeq L \hocolim L(X_\alpha). \]
\end{lemma}

\begin{proof}
It suffices to show that for any local object $Y$, there is a weak
equivalence of simplicial sets
\[ \Map^h(L (\hocolim_\alpha LX_\alpha),Y) \simeq
\Map^h(L \hocolim_\alpha X_\alpha,Y). \] This fact follows from
the following series of weak equivalences:
\[ \begin{aligned}
\Map^h (L \hocolim_\alpha LX_\alpha, Y) & \simeq
\Map^h(\hocolim_\alpha LX_\alpha,Y) \\
& \simeq \holim_\alpha \Map^h(LX_\alpha, Y) \\
& \simeq \holim_\alpha \Map^h(X_\alpha, Y) \\
& \simeq \Map^h(\hocolim_\alpha X_\alpha,Y) \\
& \simeq \Map^h(L\hocolim_\alpha X_\alpha, Y).
\end{aligned} \] \end{proof}

We now consider the simplicial space $\nerve (T_n)^t$, which is
the nerve of the free monoid on $n$ generators, considered as a
transposed simplicial space.

\begin{prop} \label{nerve}
Let $T_n$ denote the free monoid on $n$ generators.  Then in
$\LSSp_{*, c}$, $L_1 \Delta [n]^t_*$ is weakly equivalent to
$\nerve (T_n)^t$ for each $n \geq 0$.
\end{prop}

\begin{proof}
Note that when $n=0$, $\Delta [0]^t_*$ is isomorphic to $\nerve
(T_0)^t$, which is already a Segal category.

Now consider the case where $n=1$.  We want to show that the map
$\Delta [1]^t_* \rightarrow \text{nerve}(T_1)^t$ obtained by
localizing with respect to the map $\varphi_*$ is a weak
equivalence in $\LSSp_{*,c}$. In order to do this, we define a
filtration of nerve$(T_1)^t$.  Let the $k$-th stage of the
filtration be
\[ \Psi_k (\text{nerve}(T_1)^t)_j=\{(x^{n_1}|\cdots |x^{n_j})|\sum
n_j \leq k\}. \] Thus we have
\[ \Delta [1]^t_* = \Psi_1 \subseteq \Psi_2 \subseteq
\cdots \subseteq \Psi_k \subseteq \cdots \text{nerve}(T_1)^t. \]
Note that the fact that $\Delta [1]^t_* = \Psi_1$ can be observed
from looking at the bar construction notation of each as a
simplicial set (which we then view as a simplicial space by
transposing it). Each has one nondegenerate 1-simplex which we
denote by $x$.

Note that $\Psi_1$ has no nondegenerate 2-simplices.  However, we
want to define the ``composite" of $x$ with itself, a 1-simplex
which we denote by $x^2$, and add a nondegenerate 2-simplex
$[x,x]$ whose boundary consists of the 1-simplices $x, x$ and
$x^2$. More formally, if $\alpha^i$ is as defined in the
introduction, take the object
\[ G(2)^t_*= \bigcup_{i=1}^{k-1} \alpha^i \Delta[1]^t_* \subseteq \Delta[k]^t_* \]
(which lives in $\Psi_2$) and add another 1-simplex $x^2$ and the
2-simplex $[x,x]$ in $\Psi_2$. We can describe this passing from
$\Psi_1$ to $\Psi_2$ by the following pushout diagram:
\[ \xymatrix{G(2)^t_* \ar[r] \ar[d] &
\Psi_1 \ar[d] \\ \Delta [2]^t_* \ar[r] & \Psi_2} \]  Since we are
working in $\LSSp_{\ast,c}$, the left-hand vertical map is an
acyclic cofibration, and therefore $\Psi_1 \rightarrow \Psi_2$ is
an acyclic cofibration also \cite[3.14]{ds}.

Similarly, to obtain $\Psi_3$ we will add an extra 1-simplex,
denoted $x^3$, in order to add a 3-simplex $[x,x,x]$. However,
when taking the pushout, we do not want to start with $G(3)_*$,
since we have already added two of the new 1-simplices when we
localized to obtain $\Psi_2$.  So, we define $(\Delta
[3]^t_*)_{\Psi_2}$ to be the piece of $\Delta [3]^t_*$ that we
already have in $\Psi_2$. Then our pushout looks like:
\[ \xymatrix{ (\Delta [3]^t_*)_{\Psi_2} \ar[r] \ar[d] & \Psi_2 \ar[d] \\
\Delta [3]^t_* \ar[r] & \Psi_3} \] The map $(\Delta
[3]^t_*)_{\Psi_2} \rightarrow \Delta [3]^t_*$ is a weak
equivalence in $\LSSp_{\ast,c}$ as follows.  We have maps
\[ \xymatrix@1{G(3)^t_* \ar[r]^-\alpha
& (\Delta [3]^t_*)_{\Psi_2} \ar[r]^-\beta & \Delta [3]^t_*}. \]
Taking the function complex $\Map(-,X)$ for any local $X$ for any
of the three above spaces yields $X_1 \times X_1 \times X_1 \simeq
X_3$. The map $\alpha$ is a weak equivalence since it is just a
patching together of two localizations coming from the map
$G(2)^t_* \rightarrow \Delta [2]^t_*$, which is a weak equivalence
since it is one of the maps with respect to which we are
localizing. The composite map $\beta \alpha$ is also a weak
equivalence for the same reason. Thus, $\beta$ is also a weak
equivalence by the two-out-of-three property.  Again, since
$(\Delta [3]^t_\ast)_{\Psi_2} \rightarrow \Delta [3]^t_\ast$ is an
acyclic cofibration in $\LSSp_{\ast, c}$, the map $\Psi_2
\rightarrow \Psi_3$ is an acyclic cofibration also.

For greater values of $i$, define $(\Delta [i+1]^t_*)_{\Psi_i}$ to
be the piece of $\Delta [i+1]^t_*$ that we already have from
previous steps of the filtration.  Note that it is always two
copies of $\Delta [i]^t_*$ attached along a copy of $\Delta
[i-1]^t_*$, so the same argument as for $i=2$ shows that the map
$(\Delta [i+1]^t_*)_{\Psi_i} \rightarrow \Delta [i+1]^t_*$ is a
weak equivalence.  Hence, for each $i$ we obtain $\Psi_{i+1}$ via
the pushout diagram
\[ \xymatrix{(\Delta [i+1]^t_*)_{\Psi_i} \ar[r] \ar[d] & \Psi_i \ar[d] \\
\Delta [i+1]^t_* \ar[r] & \Psi_{i+1}} \]

Now that we have defined each stage of our filtration, using the
bar construction notation shows how to map our new local object to
$\text{nerve}(T_1)^t$.  For example, $[x, x, x^2] \mapsto (x, x,
x^2) \in T_1 \times T_1 \times T_1$.

Using Lemma \ref{Lhocolim} we have that
\[ \begin{aligned}
\text{nerve}(T_1)^t & \simeq L_1(\text{nerve}(T_1)^t) \\
& \simeq L_1(\hocolim (\Psi_i)) \\
& \simeq L_1(\hocolim L_1(\Psi_i)) \\
& \simeq L_1(\hocolim L_1(\Psi_1)) \\
& \simeq L_1L_1(\Psi_1)  \\
& \simeq L_1(\Psi_1) \\
& \simeq L_1(\Delta [1]^t_*).
\end{aligned} \]

Now, for $n=2$ (i.e. starting with $\Delta [2]^t_\ast$), we have
three 1-simplices, which we will call $x$, $y$, and $xy$, and one
nondegenerate 2-simplex $[x,y]$.  Because we now have two
variables, we need to define the filtration slightly differently
as $\Psi_i=\{[w_1, \ldots ,w_k]|l(w_1 \ldots w_k) \leq i \}$ where
the $w_j$'s are words in $x$ and $y$ and $l$ denotes the length of
a given word. Note that by beginning with $\Psi_1$ we start with
fewer simplices than those of the 2-simplex we are considering,
but by passing to $\Psi_2$ we obtain the $xy$ and $[x,y]$ as well
as additional nondegenerate simplices. In fact, we are actually
starting the filtration with $\Psi_1=G(2)^t_*$. The localizations
proceed as in the case where $n=1$, enabling us to map to $\nerve
(T_2)^t$.

For $n \geq 3$, the same argument works as for $n=2$, with the
filtrations being defined by the lengths of words in $n$ letters.
The resulting object is a reduced Segal category weakly equivalent
to $\Delta [n]^t_*$. Hence, we have that for any $n$, $L_1 \Delta
[n]^t_*$ is weakly equivalent to $\nerve (T_n)^t$.
\end{proof}

We now define a map $J: \Deltaop \rightarrow \Tm$ induced by the
nerve construction on a monoid $M$.  We will actually define the
map $J^{op}: {\bf \Delta} \rightarrow \Tm^{op}$.  (Note that
$\Tm^{op}$ is just the full subcategory of the category of monoids
whose objects are representatives of the isomorphism classes of
finitely generated free monoids.)

For an object $[n]$ of ${\bf \Delta}$, define $J^{op}([n]) = T_n$
where $T_n$ denotes the free monoid on $n$ generators, say $x_1,
\ldots ,x_n$.  In particular, $J^{op}([0]) = *$, the trivial
monoid.  There are coface and codegeneracy maps induced from the
nerve construction on $M$ as follows.

Taking the nerve of a simplicial monoid $M$ and regarding it as a
constant simplicial space results in a simplicial space which at
level $k$ looks like
\[ \nerve (M)_k = M^k = \Hom_{\mathcal Mon}(T_k,M). \]
By Yoneda's Lemma, the face operators
\[ d_i:\nerve(-)_k \rightarrow \nerve(-)_{k-1} \] are induced by
monoid maps $T_{k-1} \rightarrow T_k$, and similarly for the
degeneracy maps
\[ s_i:\nerve(-)_k \rightarrow \nerve(-)_{k+1}. \]  Thus the
simplicial diagram $\nerve (-)$ of representable functors $\Hom
(T_k,-)$ gives rise to a cosimplicial diagram of representing
objects $T_k$.  In particular, the coface maps are defined by:
\[ \begin{aligned} d^i(x_k) =
\begin{cases}
x_k & i<k \\
x_kx_{k+1} & i=k \\
x_{k+1} & i>k
\end{cases}
\end{aligned} \]
and the codegeneracy maps are defined analogously.  To obtain a
simplicial diagram of free monoids, we simply reverse the
direction of the arrows to obtain a functor $J: \Deltaop
\rightarrow \Tm$. This map induces a map $J^*:\SSets^{\Tm}
\rightarrow \SSets^{\Deltaop}$ which can be restricted to a map
$J^*:\SSetstm_* \rightarrow \mathcal {SS}p_*$.

We now state two definitions in the following general context. Let
$p:\mathcal C \rightarrow \mathcal D$ and $G: \mathcal C
\rightarrow \SSets$ be functors.

\begin{definition}
If $d$ is an object of $\mathcal D$, then the \emph{over category}
or \emph{category of objects over} $d$, denoted $(p \downarrow
d)$, is the category whose objects are pairs $(c,f)$ where $c$ is
an object of $\mathcal C$ and $f:p(c) \rightarrow d$ is a morphism
in $\mathcal D$.  If $c'$ is another object of $\mathcal C$, a
morphism in the over category is given by a map $c \rightarrow c'$
inducing a commutative triangle
\[ \xymatrix{p(c) \ar[dd] \ar[dr] & \\
& d \\
p(c') \ar[ur]} \]
\end{definition}

\begin{definition} \cite[11.8.1]{hirsch}
Let $p$, $c$, and $G$ be defined as above, and let $f:p(c)
\rightarrow d$ be an object in $(p \downarrow d)$. The \emph{left
Kan extension} over $p$ is a functor $p_*G: \mathcal D \rightarrow
\SSets$ defined by
\[ (p_*G)(d) = \colim_{(p \downarrow d)}((c, f)
\mapsto G(c)). \]
\end{definition}

\begin{note}
Ordinarily, we would have to take a \emph{homotopy} left Kan
extension, where we replace the colimit in the definition with a
homotopy colimit, to make sure that our calculations were homotopy
invariant.  However, since we are making our calculations in
$\LSSp_{*,c}$ where all objects are cofibrant, and since left Kan
extensions agree with homotopy left Kan extensions on cofibrant
objects \cite[3.7]{dk3}, the left Kan extension is sufficient.
\end{note}

\begin{prop} \cite[11.9.3]{hirsch}
The functor $\SSets^\mathcal C \rightarrow \SSets^\mathcal D$
given by sending $G$ the left Kan extension $p_*G$ is left adjoint
to the functor $\SSets^\mathcal D \rightarrow \SSets^\mathcal C$
given by composition with $p$.
\end{prop}


In our specific case, we define $J_*:\mathcal {SS}p_* \rightarrow
\SSets^{\Tm}_*$ to be the left Kan extension over $J$ which is, by
definition, the left adjoint to $J^*$. Note that even if $G$ is a
reduced Segal category, $J_*(G)$ is not necessarily local in
$\LSSetstm_*$. To obtain a $\Tm$-algebra, we must apply the
localization functor $L_2$.

Define $M[k]$ to be the functor $\Tm \rightarrow \SSets$ given by
$T_n \mapsto \Hom_{\Tm}(T_k,T_n)=T_k^n$.
Let $G$ be the reduced Segal category $\nerve
(T_k)^t$.

\begin{lemma} \label{mk}
In $\LSSetstm_*$, $L_2J_*(G)$ is weakly equivalent to $M[k]$.
\end{lemma}

\begin{proof}
It suffices to show that for any local object $X$ in
$\LSSetstm_*$,
\[\Map^h_{\LSSetstm_*}(L_2J_*G,X) \simeq X(T_k) \] since $\Map^h_{\LSSetstm}(M[k],X)$ is precisely $X(T_k)$.
This fact can be shown in the following argument:
\[ \begin{aligned}
\Map^h_{\LSSetstm_*}(L_2J_*G,X)& \simeq \Map^h_{\LSSetstm_*}(J_*G,X) \\
& \simeq \Map^h_{\mathcal {LSS}p_{*,c}}(G,J^*X) \\
& \simeq \Map^h_{\mathcal {LSS}p_{*,c}}(L_1 \Delta [k]^t_*, J^*X) \\
& \simeq \Map^h_{\mathcal {LSS}p_{*,c}}(\Delta [k]^t_*, J^*X) \\
& \simeq J^*X[k] \\
& \simeq X(T_k).
\end{aligned} \]
\end{proof}


\begin{prop} \label{JLJ}
For any object $X$ in $\mathcal{SS}p_{*,c}$, we have that $L_1X$
is weakly equivalent to $J^*L_2J_*X$.
\end{prop}

\begin{proof}
First note that $X \simeq \hocolim_{\Deltaop}([n] \rightarrow
\coprod_i \Delta [n_i]^t)$ where the values of $i$ depend on $n$.
We begin by looking at $L_1 X$.  Using Lemma \ref{Lhocolim}, we
have the following:
\[ \begin{aligned}
L_1X & \simeq L_1 \hocolim_{\Deltaop} ([n] \mapsto \coprod \Delta
[n_i]^t_*)
\\
& \simeq L_1 \hocolim_{\Deltaop} L_1([n] \mapsto \coprod \Delta [n_i]^t_*) \\
& \simeq L_1 \hocolim_{\Deltaop} ([n] \mapsto \text{nerve}(T_{\sum
n_i})^t)
\end{aligned} \]
However, $\hocolim_{\Deltaop} ([n] \mapsto \nerve(T_{\sum
n_i})^t)$ is already local, a fact which follows from the fact
that the homotopy colimit can be taken at each level, yielding a
Segal precategory in $\mathcal{LSS}p_{*,c}$ which is still a Segal
category.

Working from the other side of the desired equation, we obtain,
using the fact that left adjoints commute with colimits and
homotopy colimits:
\[ \begin{aligned}
J^*L_2J_*X & \simeq J^*L_2J_* \hocolim_{\Deltaop}([n] \mapsto
\coprod
\Delta [n_i]^t_*) \\
& \simeq J^*L_2 \hocolim_{\Deltaop} J_*([n] \mapsto \coprod \Delta [n_i]^t_*) \\
& \simeq J^*L_2 \hocolim_{\Deltaop} L_2J_*([n] \mapsto \coprod \Delta [n_i]^t_*) \\
& \simeq J^*L_2 \hocolim_{\Deltaop} ([n] \mapsto M[\sum n_i]) \\
\end{aligned} \]

At each level, we have the same spaces as in the
$\Deltaop$-diagram, but with maps given by $\Tm$ rather than
$\Deltaop$.  Thus, applying the restriction map $J^*$ results in a
diagram with the same objects at each level, as we wished to show.
\end{proof}

We begin by showing that we have a Quillen pair even before we
apply the localization functors.  Notice that we return to using
our model structure $\SSp_{*,f}$ with levelwise fibrations.

\begin{prop}
The adjoint pair
\[ \xymatrix@1{J_*: \mathcal {SS}p_{*,f} \ar@<.5ex>[r] & \SSets^{\Tm}_* :J^* \ar@<.5ex>[l]} \]
is a Quillen pair.
\end{prop}

\begin{proof}
In each case the fibrations and weak equivalences are defined
levelwise.  Since the right adjoint $J^*$ preserves the spaces at
each level, it must preserve both fibrations and acyclic
fibrations.
\end{proof}

We then need to show that this Quillen pair induces a Quillen
equivalence between the localized model category structures in
order to prove our main theorem.  While we will be working with
$\LSSp_{*,f}$ rather than $\LSSp_{*,c}$, the calculations we made
above will still hold because the weak equivalences are the same
in each case (see the end of section 3).

\begin{proof}[Proof of Theorem \ref{main}]
We first need to show that we still have a Quillen pair even after
localizing each category.  This fact follows from Theorem
\ref{LocPair} after we notice that the left derived localizing set
${\bf L}J_* \varphi$ \cite[8.5.11]{hirsch} is the same as the set
of maps $S=\{p_0, p_1, \ldots \}$ that we localize with respect to
in order to obtain $\LSSetstm_*$, as follows.  Consider the maps
$\varphi_k: \coprod_k \Delta [1]^t_\ast \rightarrow \Delta
[k]^t_\ast$ for each $k$.  If we localize this map, we obtain
$\coprod_k \nerve (T_1)^t \rightarrow \nerve (T_k)^t$, where $T_k$
denotes the free monoid on $k$ generators.  Then we can apply the
functor $J_\ast$ to obtain a map
\[ J_\ast(\coprod_k \nerve (T_1)^t) \rightarrow J_\ast (\nerve
(T_k)^t) \] which we can then localize to obtain, by Lemma
\ref{mk}, the map $\coprod_k M[1] \rightarrow M[k]$.  However,
this is precisely the map $p_k:\coprod_k \Hom (T_1,-) \rightarrow
\Hom (T_k,-)$.

We now need to show that we have a Quillen equivalence.  First, we
need to know that the right adjoint $J^*$ reflects weak
equivalences between fibrant objects.  In each of the two
localized model categories $\LSSp_{\ast, f}$ and $\LSSetstm_*$, an
object is fibrant if and only if it is local and fibrant in the
unlocalized model category. Therefore, in each case a weak
equivalence between fibrant objects is a levelwise weak
equivalence.  Since $J^*$ does not change the spaces at each
level, it must reflect weak equivalences between fibrant objects.

Finally, by Proposition \ref{JLJ}, $L_1X \simeq J^*L_2J_*X$ for
any functor $X: \Deltaop \rightarrow \SSets$, and in particular
for any cofibrant $X$.
\end{proof}

Notice that our Quillen equivalences compose:
\[ \LSSp_{*,f} \rightleftarrows \LSSetstm_* \rightleftarrows \LSSetstm \rightleftarrows
\Algtm \] (where the left adjoint functors are the topmost maps).
Therefore, we actually have a single Quillen equivalence
\[ \LSSp_{*,f} \rightleftarrows \Algtm. \]

\section{A Generalization to Simplicial Categories and Segal
Categories}

We would like to generalize our result to the category of
simplicial categories which have a set $\mathcal O$ of objects,
where the morphisms are required to be the identity on this set,
since a simplicial category is a generalization of a simplicial
monoid.

\begin{definition}
A \emph{simplicial category} is a category enriched over
simplicial sets, i.e. a category in which the morphisms between
any two objects form a simplicial set.
\end{definition}

Note that a simplicial category in our sense is a simplicial
object in the category of categories in which the categories at
each level have the same objects.  In particular, the face and
degeneracy maps are the identity on the objects.

Analogously, we will generalize to Segal categories with the same
set $\mathcal O$ in degree zero.  As we have seen, any Segal
category can be described as a functor from the simplicial
category $\Deltaop$ to the category $\SSets$ of simplicial sets
satisfying two conditions: discreteness of the zero space and a
product condition up to homotopy (Definition \ref{SeCat}).
However, we will find it more convenient, when dealing with Segal
categories with a fixed object set, to use a larger indexing
category which specifies the objects.

Let $\mathcal O$ be a set.  We define the category
$\Deltaop_{\mathcal O}$ as follows.  The objects are given by
$[n]_{x_0, \ldots ,x_n}$ where $n \geq 0$ and $(x_0, \ldots x_n)
\in \mathcal O^{n+1}$.  The $[n]$ should be thought of as in the
simplicial category $\Deltaop$; however, recall that when we work
with Segal categories we will require all morphisms to preserve
the objects. Therefore, we need to have a separate $[n]$ for each
possible $(n+1)$-tuple of objects in $\mathcal O$. The morphisms
in $\Deltaop_\mathcal O$ are those of $\Deltaop$ but depend on the
choice of $(x_0, \ldots ,x_n)$.  Specifically, the face maps are
\[ d_i:[n]_{x_0, \ldots ,x_n} \rightarrow [n-1]_{x_0, \ldots
,\widehat x_i, \ldots ,x_n} \] and the degeneracy maps are
\[ s_i:[n]_{x_0, \ldots ,x_n} \rightarrow [n+1]_{x_0, \ldots,
x_{i-1}, x_i, x_i, x_{i+1}, \ldots ,x_n}. \]  We will sometimes
refer to our Segal precategories with $\mathcal O$ in degree zero
as $\Deltaop_\mathcal O$-spaces.  Note that if $\mathcal O$ is the
one-object set, then $\Deltaop_{\mathcal O}$ is just $\Deltaop$.

Now we can use this notation to describe Segal categories.  So, a
Segal category with $\mathcal O$ in degree zero is a functor
$X:\Deltaop_\mathcal O \rightarrow \SSets$ such that for each $k
\geq 2$, the map
\[ X([n]_{x_0, \ldots, x_n}) \rightarrow X([1]_{x_0, x_1}) \times
_{X[0]_{x_1}} \cdots \times_{X[0]_{x_{n-1}}} X([1]_{x_{n-1},x_n})
\] is a weak equivalence.

Recall from section 3 that we have model category structures
$\LSSpof$ and $\LSSpoc$ on the category of Segal precategories
with $\mathcal O$ in degree zero, in each of which the fibrant
objects are Segal categories.

We would like to think of the category of simplicial categories
with object set $\mathcal O$ as a diagram category as well.  To do
so, we need to define the notion of a multi-sorted algebraic
theory. To see more details, see \cite{multisort}.

\begin{definition}
Given a set $S$, an $S$-\emph{sorted algebraic theory} (or
\emph{multi-sorted theory}) $\mathcal T$ is a small category with
objects $T_{\alphau^n}$ where $\alphau^n = <\alpha_1, \ldots
,\alpha_n>$ for $\alpha_i \in S$ and $n \geq 0$ varying, and such
that each $T_{\alphau^n}$ is equipped with an isomorphism
\[ T_{\alphau^n} \cong \prod_{i=1}^n T_{\alpha_i}. \]
For a particular $\alphau^n$, the entries $\alpha_i$ can repeat,
but they are not ordered.  There exists a terminal object $T_0$
(corresponding to the empty subset of $S$).
\end{definition}

In particular, we can talk about the theory of $\mathcal
O$-categories, which we will denote by $\Tocat$.  To define this
theory, first consider the category $\mathcal {OC}at$ whose
objects are the categories with a fixed object set $\mathcal O$
and whose morphisms are the functors which are the identity map on
the objects. The objects of $\Tocat$ are categories which are
freely generated by directed graphs with vertices corresponding to
the elements of the set $\mathcal O$. This theory will be sorted
by pairs of elements in $\mathcal O$, corresponding to the
morphisms with source the first element and target the second.  In
other words, this theory is $(\mathcal O \times \mathcal
O)$-sorted \cite[3.5]{multisort}. (In the one-object case, we have
the ordinary theory of monoids.) We can then say that a simplicial
category with object set $\mathcal O$ is a strict
$\Tocat$-algebra, where strict and homotopy $\mathcal T$-algebras
are defined for multi-sorted theories $\mathcal T$ just as for
ordinary algebraic theories.

Again, we have a model structure $\Algt$ on the category of all
$\mathcal T$-algebras and a model structure $\SSetst$ on the
category of all functors $\mathcal T \rightarrow \SSets$ which can
be localized as before to obtain a model category structure
$\LSSetst$ in which the local objects are homotopy $\mathcal
T$-algebras \cite{multisort}.  For $\Tocat$, we can define a
category $\SSetsTocat_\mathcal O$ of functors $\Tocat \rightarrow
\SSets$ which send $T_0$ to $\coprod_\mathcal O \Delta [0]$.
Making modifications as in the case of $\LSSetstm_*$, we can
define a model structure $\LSSetsTocat_\mathcal O$ which is
Quillen equivalent to $\LSSetsTocat$.

Furthermore, Theorem \ref{rigid} for algebraic theories can be
generalized to the case of multi-sorted theories, and therefore we
have that there is a Quillen equivalence of model categories
between $\Algt$ and $\LSSetst$ for any multi-sorted theory
$\mathcal T$ \cite{multisort}. Applying this result to the
multi-sorted theory $\Tocat$ and using the Quillen equivalence in
the previous paragraph, we have reduced the problem to finding a
Quillen equivalence between $\LSSpof$ and
$\mathcal{LSS}ets^{\Tocat}_\mathcal O$.  We begin as before by
defining a map of diagrams $J:\Deltaop_{\mathcal O} \rightarrow
\Tocat$ which then induces a map $J^*:\SSetsTocat_\mathcal O
\rightarrow \SSetso$.


Given an integer $n \geq 0$ and $\xu = (x_0, \ldots ,x_n) \in
\mathcal O^{n+1}$, define $T_{n, \xu}$ to be the free category
with object set $\mathcal O$ and morphisms freely generated by the
set $\{x_{i-1} \rightarrow x_i | 1 \leq i \leq n\}$. However, note
that this free category is, in light of the definition of
multi-sorted theory, just the object $T_\alphau$ where $\alphau =
(\alpha_1, \ldots ,\alpha_n)$ and $\alpha_i = (x_{i-1},x_i) \in
\mathcal O \times \mathcal O$.  (We will find both notations
useful.)

We then obtain a simplicial diagram as follows:
\[ \Hom(T_0, \mathcal C) \Leftarrow \coprod_{\alpha_1}
\Hom(T_{\alpha_1} ,\mathcal C) \Lleftarrow \coprod_{\alphau =
(\alpha_1,\alpha_2)} \Hom(T_\alphau ,\mathcal C) \cdots \] Because
we have a coproduct at each level, this diagram is not
representable by a cosimplicial diagram.  However, if we split up
the coproducts and consider it as a diagram of objects
$\Hom(T_\alphau ,\mathcal C)$ varying over $\alphau$ and $n$, we
get a map $J:\Deltaop_\mathcal O \rightarrow \Tocat$. Then we
define the map
\[ J^*: \SSetsTocat_\mathcal O \rightarrow \SSetso \]
to be the one induced by this map $J$ above.  Using the model
structures on these categories we get
\[ J^*: \LSSetsTocat_\mathcal O \rightarrow \LSSpof. \]

We will denote by $L_1$ the localization of simplicial spaces with
$\mathcal O$ in dimension zero, namely the functorial fibrant
replacement functor in $\LSSpoc$.  (As above, we will make our
calculations in this category rather than in $\LSSpof$.)  We have
an analogous localization in the category $\LSSetsTocat_\mathcal
O$, which we will denote by $L_2$.

From there, we can apply the methods used in the previous section.
We first consider what happens when we localize an $n$-simplex
$\Delta [n]_{x_0, \ldots ,x_n}$. If $x_i = x_j$ for all $0 \leq
i,j \leq n$ then we get the free monoid on $n$ generators, with
object $x_i$. However, if $x_i \neq x_j$ for some $0 \leq i,j \leq
n$, we get the nerve of the free category generated by the
directed graph
\[ x_0 \rightarrow x_1 \rightarrow \cdots \rightarrow x_n. \]
For example, if $x_i \neq x_j$ for all $0 \leq i,j \leq n$, then
the localization only consists of including in composites, for
example the 1-simplices $x_i \rightarrow x_{i+2}$.  If $n=1$ and
$x_0 \neq x_1$, then the localization does not change the Segal
category since there are no ``composable" 1-simplices.  We will
denote the free category by $T_{n,\xu}$, as before, and its nerve
by $G_{n,\xu} = \nerve(T_{n,\xu})^t$.

Now, define $J_*$ be the left Kan extension of the map $J^*$.
However, note that again we must apply the localization functor
$L_2$ to assure that we get a homotopy $\Tocat$-algebra in
$\LSSetsTocat_\mathcal O$. As before, it suffices to show that for
any $\Deltaop_\mathcal O$-space $X$, we have that $L_1X =
J^*L_2J_*X$. We will first need a lemma.

Analogous to $M[k]$ above, define $C[n]_\xu$ to be the
$\Tocat$-diagram in of simplicial sets which has the same
simplicial set at each level as $G_{n, \xu}$, but with theory maps
between them rather than the simplicial maps.

\begin{lemma}
In $\LSSetsTocat_\mathcal O$, $L_2J_*(G_{n,\xu})$ is weakly
equivalent to $C[n]_\xu$.
\end{lemma}

\begin{proof} It suffices to show that for a local
object $X$ in $\LSSetsTocat_\mathcal O$, we have that
\[ \Hom_{\LSSetsTocat_\mathcal O}
(L_2J_*(G_{n,\xu}),X) \simeq X(T_{n,\xu}). \] This can be proved
as follows:
\[ \begin{aligned}
\Hom_{\LSSetsTocat_\mathcal O} (L_2J_*(G_{n,\xu}),X) & \simeq
\Hom_{\LSSetsTocat_\mathcal O}(J_*(G_{n,\xu}),X) \\
& \simeq \Hom_{\LSSpoc}(G_{n,\xu}, J^*X) \\
& \simeq \Hom_{\LSSpoc}(L_1 \Delta[n]_\xu, J^*X) \\
& \simeq \Hom_{\LSSpoc}(\Delta[n]_\xu, J^*X) \\
& \simeq J^*(X)[n]_\xu \\
& \simeq X(T_{n,\xu})
\end{aligned} \]
\end{proof}

Now, using this lemma, the following proposition suffices to prove
the theorem.

\begin{prop}
Given any $\Deltaop_\mathcal O$-space $X$, $L_1X$ and $J^*L_2J_*X$
are weakly equivalent.
\end{prop}

\begin{proof}
Given any $\Deltaop_\mathcal O$-space $X$, it can be written as
the homotopy colimit
\[ X \simeq \hocolim_{\Deltaop}([n] \mapsto \amalg_{n_i,\xu} \Delta
[n_i]_\xu)\] where the $n_i$ depend on $n$ and $\xu = (x_1, \ldots
,x_{n_i})$.  We then have the following:
\[ \begin{aligned}
L_1X & \simeq L_1 \hocolim_{\Deltaop}([n] \mapsto \amalg_{n_i,\xu}
\Delta [n_i]_\xu) \\
& \simeq L_1 \hocolim_{\Deltaop} L_1([n] \mapsto \amalg_{n_i,\xu}
\Delta [n_i]_\xu) \\
& \simeq L_1 \hocolim_{\Deltaop} ([n] \mapsto \amalg_{n_i,\xu}
G_{n_i,\xu}).
\end{aligned} \]
Working from the other side, we have
\[ \begin{aligned}
J^*L_2J_*X & \simeq J^*L_2J_* \hocolim_{\Deltaop}([n] \mapsto
\amalg_{n_i,\xu} \Delta [n_i]_\xu) \\
& \simeq J^*L_2 \hocolim_{\Deltaop} J_*([n] \mapsto
\amalg_{n_i,\xu}
\Delta [n_i]_\xu) \\
& \simeq J^*L_2 \hocolim_{\Deltaop} L_2J_*([n] \mapsto
\amalg_{n_i,\xu}
\Delta [n_i]_\xu) \\
& \simeq J^*L_2 \hocolim_{\Deltaop}([n] \mapsto \amalg_{n_i,\xu}
C[n_i]_\xu).
\end{aligned} \]
But these two are equal, by the above lemma.
\end{proof}

Now, the following theorem can be proved just as Theorem
\ref{main} was in the previous section.

\begin{theorem} \label{cats}
The adjoint pair
\[ \xymatrix@1{J_*: \LSSpof \ar@<.5ex>[r] & \LSSetsTocat_\mathcal O :J^* \ar@<.5ex>[l]} \]
is a Quillen equivalence.
\end{theorem}

In particular, composing with the Quillen equivalences given by
the generalizations of Theorem \ref{rigid} and Proposition
\ref{tmmcs}, there is a Quillen equivalence
\[ \LSSpof \rightleftarrows \mathcal Alg^{\Tocat}. \]

\bibliographystyle{amsalpha}

\end{document}